\documentstyle[12pt]{article}
\begin{document}
\large
\newtheorem{lem}{Lemma}
\newtheorem{prop}{Proposition}
\newtheorem{cor}{Corollary}
\newtheorem{teo}{Theorem}
\title{Limit Theorems for Sums of $p$-Adic Random Variables}
\author{Anatoly N. Kochubei\\ \footnotesize Institute of Mathematics,\\
\footnotesize Ukrainian National Academy of Sciences,\\
\footnotesize Tereshchenkivska 3, Kiev, 252601 Ukraine}
\date{}
\maketitle
\begin{abstract}
We study $p$-adic counterparts of stable distributions, that is
limit distributions for sequences of normalized sums of
independent identically distributed $p$-adic-valued random
variables. In contrast to the classical case, non-degenerate
limit distributions can be obtained only under certain
assumptions on the asymptotic behaviour of the number of
summands in the approximating sums. This asymptotics determines
the ``exponent of stability''.
\end{abstract}
\section{Introduction}

The studies in infinitesimal systems of probability measures on
locally compact groups (see [He] and references therein) are
concentrated mainly on the problem of convergence to the
Gaussian distribution. Not much is known (see a review in [Kh])
about analogues of
stable distributions and the systems converging to them. Of
course this is connected with the fact that standard
normalization procedures (see e.g. [F]) are not possible for
general groups.

This paper is devoted to an important example of a group for
which an analogue of the classical theory can be constructed
though the results are quite different from the ones for the
group ${\bf R}$. Namely, we shall consider the additive group
of the field $Q_p$ of $p$-adic numbers.

Note that there is no Gaussian measure on $Q_p$ (in the sense
of Parthasarathy) since $Q_p$ is totally disconnected. On the
other hand, the distributions $G_{a,\alpha }$  on $Q_p$ having
the functions $g_{a,\alpha }(t)=\exp (-a|t|_p^\alpha ) $ ,
$a>0,\ \alpha >0$ ($|\cdot |_p$ is the $p$-adic absolute value;
see Sect. 2) as their Fourier transforms, were used recently by
several authors ([B], [I], [Ha], [K1], [K2], [Va], [VVZ]) as $p$-adic
counterparts of stable distributions, being the basis of the
$p$-adic stochastic analysis initiated in the above papers and
related to $p$-adic models of mathematical physics; for other
approaches to $p$-adic stochastic processes see [AK], [E1], [E2].

It is natural to try to obtain these and more general distributions
``of stable type'' as limits of certain normalized sums of
independent identically distributed random variables with values
in $Q_p$. Even the above model example shows differences
between the  $p$-adic and real cases. Since the  $p$-adic
absolute value can equal only an integer power of $p$, the
classical definitions of a stable distribution (as well as the
generalizations proposed in [T]) do not make sense for $Q_p$.

Let $X_1,\ldots ,X_n,\ldots $ be a sequence of
independent identically distributed  $Q_p$-valued random
variables, and $B_1,\ldots ,B_n,\ldots $ a sequence of $p$-adic
numbers. We consider the normalized sums
\begin{equation}
S_n=B_n^{-1}\left(X_1+\cdots +X_{k(n)}\right) \ ,\ \ n=1,2,\ldots
\end{equation}
where $\{k(n)\}$ is an increasing sequence of natural numbers.
Let $F_n$ be a distribution of $S_n$; suppose that $F_n\to G$
in the weak sense. Our main aim is to describe some distributions
$G$ (or their characteristic functions $g(t)$ ) which may
appear this way. Thus we confine ourselves to the ``strictly
stable'' case. A significance of centering is not clear for
$p$-adic random variables (for which, by the way, an expectation
is not defined).

If $|B_n|_p\le \mbox{const}$ then the distribution $G$ is either
degenerate (equal to the delta measure) or equal to a cutoff of
the Haar measure of the additive group of $Q_p$. An interesting
case is the one when $|B_n|_p\to \infty $ , so that $G$ is
infinitely divisible [PRV]. The answer depends on the behaviour
of the sequence
$$
\rho _n=\frac{k(n)}{k(n+1)}\ ,\ \ n=1,2,\ldots .
$$
Passing to subsequences we may assume that $\rho _n\to \beta $
, $0\le \beta \le 1$.

If $\beta =1$ (as in the 'classical" case $k(n)=n$) then, in
sharp contrast to the case of real-valued random variables, $G$
is degenerate. Another extreme case is $\beta =0$ when either
$G$ is degenerate or $g$ has a compact support.

The $p$-adic counterparts of stable distributions emerge when
$0<\beta <1$. These include $G_{a,\alpha }$ (for which $\beta
=p^{-\alpha }$ ). We find a class of distributions (defined by
a functional equation for their L\'evy measures) which
correspond to weak
limits of sequences (1). Its subclass consisting of symmetric
distributions coincides with the set of distributions
corresponding to weak limits of sequences (1) for which the
random variables appearing in (1) have symmetric distributions.
We have not found such a complete description for the
non-symmetric case. The difficulty here may be related to the
fact of non-uniqueness of the L\'evy-Khinchin representation of
an infinitely divisible distribution on $Q_p$ (as on any
Abelian group possessing compact subgroups; see [PRV]).

Finally, we describe the domains of attraction for the above
distributions giving conditions for the weak convergence of the
sequence (1).

The author is grateful to E.D.Belokolos who called the author's
attention to the problem, for useful discussions.

\section{Preliminaries}

In this section we give some basic information from $p$-adic
analysis. See [VVZ] for further details.

Let $p$ be a prime number. The field of $p$-adic numbers is the
completion $Q_p$ of the field of rational numbers, with respect
to the absolute value $|x|_p$ defined by setting $|0|_p=0$ and
$|x|_p=p^{-\nu }$ if $x=p^\nu \displaystyle\frac{m}{n}$ where $\nu ,m,n\in
{\bf Z}$ and $m,n$ are prime to $p$.

The absolute value $|x|_p\ ,\ x\in Q_p$, has the following
properties: $|x|_p=0$ if and only if $x=0$; $|xy|_p=|x|_p|y|_p$;
$|x+y|_p\le \max (|x|_p\ ,\ |y|_p)$. If $|x|_p=p^N$ then $x$
admits the canonical representation
\begin{equation}
x=p^{-N}\left( x_0+x_1p+x_2p^2+\cdots \right)
\end{equation}
where $x_0,x_1,x_2,\ldots \in \{0,1,\ldots ,p-1\},\ x_0\ne 0$.
The series is convergent with respect to the topology defined
by the metric $|x-y|_p$.

$Q_p$ is a complete, separable, totally disconnected, locally
compact metric space. We shall denote by $dx$ the Haar measure
on the additive group of $Q_p$ normalized in such a way that
$$
\int \limits _{|x|_p\le 1}dx=1.
$$
If $a\in Q_p\ ,\ a\ne 0$, then $d(xa)=|a|_p\,dx$. The measure
of a ball $\{x\in Q_p\ :\ |x|_p\le p^N\}$ equals $p^N$. Note
that a ball, as well as a sphere $\{x\in Q_p\ :\ |x|_p= p^N\}$,
are open and simultaneously closed (compact) sets.

The canonical additive character of the field $Q_p$ is defined
by the formula
$$
\chi (x)=\exp \left( 2\pi i\{x\}_p\right)
$$
where $\{x\}_p$ is the fractional part of $x\in Q_p$; if $x$
has the representation (2) then
$$
\{x\}_p=p^{-N}\left( x_0+x_1p+\cdots +x_{N-1}p^{N-1}\right)
$$
if $N>0$, and $\{x\}_p=0$ if $N\le 0$.

The character $\chi $ is an example of a locally constant
function on $Q_p$ : $\chi (x+x')=\chi (x)$ for any $x\in Q_p$,
if $|x'|_p\le 1$. In general a function $f: Q_p\to {\bf C}$
is called locally constant if there exists such $n\in {\bf Z}$
that $f(x+x')=f(x)$ for any $x\in Q_p$, if $|x'|_p\le p^n$.

The Fourier transform of a complex-valued function $\varphi\in
L_1(Q_p)$ is defined by
\begin{equation}
\widehat{\varphi }(\xi )=\int \limits _{Q_p}\chi (\xi x)\varphi (x)
\, dx\ ,\ \ \xi \in Q_p\ .
\end{equation}
The inverse transform is
\begin{equation}
\varphi (x)=\int \limits _{Q_p}\chi (-\xi x)\widehat{\varphi }(\xi )
\, d\xi \ ,\ \ x\in Q_p\ ,
\end{equation}
if $\widehat{\varphi }\in L_1(Q_p)$. In particular, the relations
(3), (4) are valid for $\varphi \in {\cal D}(Q_p)$ where
${\cal D}(Q_p)$ is the space of locally constant functions with
compact supports. In this case $\varphi \in {\cal D}(Q_p)$
implies $\widehat{\varphi }\in {\cal D}(Q_p)$. Note that
${\cal D}(Q_p)$ contains, in particular, indicator functions of
all open compact subsets of $Q_p$.

Let $\mu $ be a probability measure on the Borel $\sigma
$-algebra of $Q_p$. Its characteristic function is defined as
usual:
$$
\widehat{\mu }(t)=\int \limits _{Q_p}\chi (tx)\mu (dx).
$$
If $\mu $ is symmetric, that is $\mu (-M)=\mu (M)$ for any
Borel set $M\subset Q_p$, then $\widehat{\mu }$ is a real-valued
function. If $\mu =\delta _\xi $ is a delta measure
concentrated at $\xi \in Q_p$ then $\widehat{\mu }(t)=\chi (\xi
t)$.

\begin{lem}
If $|\widehat{\mu }(t_0)|=1$ for some $t_0\in Q_p,\ t_0\ne 0$,
then $\widehat{\mu }$ is a
locally constant function. If $|\widehat{\mu }(t)|$ takes only two
values, 0 and 1, then there exists $\xi \in Q_p$,\linebreak $N\in {\bf
Z}$ such that
\begin{equation}
\widehat{\mu }(t)=\chi (t\xi )\Omega _N(t)
\end{equation}
where
$$
\Omega _N(t)=\left\{ \begin{array}{rl}
1, & \mbox{if }\ \ |t|_p\le p^N \\
0, & \mbox{if }\ \ |t|_p>p^N
\end{array} \right.
$$
In this case $\mu (dx)=p^N\Omega _{-N}(x-\xi )dx$. If $|\widehat{\mu
}(t)|\equiv 1$ then $\widehat{\mu }(t)=\chi (t\xi ),\ \mu =\delta
_\xi\ ,\ \xi \in Q_p$.
\end{lem}

{\it Proof.} Let $t_0\in Q_p$ be such that $t_0\ne 0,\
|\widehat{\mu }(t_0)|=1$. Denote by $R_p$ the set of rational
numbers of the form $p^{-n}\left( a_0+a_1p+\cdots
+a_{n-1}p^{n-1}\right)$, $n\ge 1$; $a_0,\ldots ,a_{n-1}\in
\{0,1,\ldots ,p-1\}$, $a_0\ne 0$.

Suppose that $\widehat{\mu }(t_0)=e^{2\pi ir},\ 0\le r<1$. Then by
the definition of $\chi $
$$
e^{2\pi ir}=\int \limits _{Q_p}\exp (2\pi i\{t_0x\}_p)\mu (dx)
$$
whence
$$
\int \limits _{Q_p}(1-\exp (2\pi i\{t_0x\}_p-r)))\mu (dx)=0.
$$
In particular,
$$
\int \limits _{Q_p}(1-\cos 2\pi (\{t_0x\}_p-r))\mu (dx)=0.
$$
This means that either $r\in R_p$ or $r=0$.

In both cases there exists $\xi \in Q_p$ such that
$r=\{t_0\xi\}_p$ , that is $\widehat{\mu }(t_0)=\chi (t_0\xi )$. As
above, we obtain that
$$
\int \limits _{Q_p}(1-\cos 2\pi (\{t_0(x-\xi )\}_p)\mu (dx)=0,
$$
so that $\mu $ is concentrated on the set of those $x$ for
which $\{t_0(x-\xi )\}_p=0$, that is on the set
$\left\{x\in Q_p:\ \ |x-\xi |_p\le |t_0|_p^{-1}\right\}$.

Now
$$
\widehat{\mu }(t)=\int \limits _{|x-\xi |_p\le |t_0|_p^{-1}}\chi
(tx)\mu (dx)=\chi (t\xi )\int \limits _{|x-\xi |_p\le
|t_0|_p^{-1}}\chi (t(x-\xi ))\mu (dx)
$$
so that $\widehat{\mu }$ is locally constant and $\widehat{\mu
}(t)=\chi (t\xi )$ if $|t|_p\le |t_0|_p$.

Suppose that  $|\widehat{\mu }(t)|=0$ or 1. If the set $\{t\in Q_p\
: \ |\widehat{\mu }(t)|=1\}$ is unbounded then it coincides with
$Q_p$ , and in this case $\widehat{\mu }(t)=\chi (t\xi ),\ \mu =\delta
_\xi $. Otherwise we come to (5). The expression for $\mu $
follows from well-known integration formulas (see [VVZ]). $\Box
$

\vspace{\baselineskip}
If $\mu $ is an infinitely divisible distribution then it
follows from the general result of [PRV] that
\begin{equation}
\widehat{\mu }(t)=\chi (t\xi )\Omega _N(t)\exp \int \limits
_{Q_p\setminus \{0\}}(\chi (tx)-1)\Phi (dx),\ \ \ t\in Q_p\ ,
\end{equation}
where $\xi \in Q_p\ ,\ N\in {\bf Z}\cup \{\infty \}\ \ (\Omega
_\infty (t)\equiv 1)$, $\Phi $ is a Borel measure on
$Q_p\setminus \{0\}$ which is finite on the complement of any
neighbourhood of zero. Formula (6) differs from a similar
formula for ${\bf R}$ in two respects - possible presence of
the factor $\Omega _N(t)$ (thus $\widehat{\mu }$ may vanish on an
open set), and non-uniqueness of the L\'evy measure $\Phi $.
However, $\Phi $ can be uniquely determined if the integral
under the exp is given.

\begin{lem}
If
\begin{equation}
\varphi (t)=\int \limits
_{Q_p\setminus \{0\}}(\chi (tx)-1)\Phi (dx),\ \ \ t\in Q_p\ ,
\end{equation}
then for any open compact subset $M\subset Q_p\setminus \{0\}$
\begin{equation}
\Phi (M)=\int \limits _{Q_p}\varphi (y)m(y)\,dy
\end{equation}
where $m$ is an inverse Fourier transform of the indicator
function $\omega _M$ of the set $M$.
\end{lem}

{\it Proof.} We have
$$
\omega _M(x)=\int \limits _{Q_p}\chi (xy)m(y)\,dy\ ,\ \ \ x\in
Q_p\ ,
$$
whence $m\in {\cal D}(Q_p)$ and
\begin{equation}
\int \limits _{Q_p}m(y)\,dy=\omega _M(0)=0.
\end{equation}

Using (7) and (9) we obtain that
$$
\int \limits _{Q_p\setminus \{0\}}\omega _M(x)\Phi (dx)= \int \limits
_{Q_p\setminus \{0\}}\Phi (dx)\int \limits _{Q_p}(\chi (xy)-1)m(y)\,dy
=\int \limits _{Q_p}\varphi (y)m(y)\,dy
$$
which is equivalent to (8). Our use of the Fubini theorem was
based on the fact that supp $m\subset \{y\in Q_p\ :\ |y|_p\le
p^l\}$ for some $l\in {\bf Z}$, and $\chi (xy)-1=0$ for
$|y|_p\le p^l,\ |x|_p\le p^{-l}$ while the measure $\Phi $ is
finite on the set $\{x\in Q_p\ :\ |x|_p>p^{-l}\}$. $\Box $

\section{Limits of Normalized Sums}

Let us consider the normalized sums (1) with $|B_n|_p\to \infty
$ and $\rho _n\to \beta \ ,\ 0\le \beta \le 1$. Let $\gamma
_n=\frac{B_n}{B_{n+1}}$. Since $|B_n|_p\to \infty $, there
exists a subsequence $\{\gamma _{n_l}\}$ for which
$|\gamma _{n_l}|_p\le p^{-1}$. We may assume (passing if
necessary to a subsequence once more) that $\gamma _n\to \gamma
_0$ in $Q_p$, $|\gamma _0|_p\le p^{-1}$.

Suppose that the distributions $F_n$ of the normalized sums
$S_n$ converge weakly, $F_n\to G$, and $g(t)$ is a
characteristic function of $G$. Let $f(t)$ be the
characteristic function of each of the (independent,
identically distributed) random variables $X_n$. Then
\begin{equation}
\left( f\left( \frac{t}{B_n}\right) \right) ^{k(n)}\to g(t)\ ,\
\ n\to \infty \ ,
\end{equation}
uniformly on compact subsets of $Q_p$. The left-hand side of
(10) will be denoted $f_n(t)$.

\begin{prop}
(i) If $\beta \ne 0$ then uniformly on compact subsets of $Q_p$
\begin{equation}
|f_n(\gamma _nt)|\to |g(t)|^\beta .
\end{equation}
If $\beta =0$, this relation holds for those $t$ where $g(t)\ne
0$.\\
(ii) The identity
\begin{equation}
|g(\gamma _0t)|=|g(t)|^\beta
\end{equation}
is valid for any $t\in Q_p$, if $\beta \ne 0$, and for any $t$
with $g(t)\ne 0$ if $\beta =0$.
\end{prop}

{\it Proof.} Let us consider a random variable
$S_n'=B_{n+1}^{-1}\left( X_1+\cdots +X_{k(n)}\right) $. Its
characteristic function equals
\begin{equation}
\left( f\left( \frac{t}{B_{n+1}}\right) \right)
^{k(n)}=f_n(\gamma _nt).
\end{equation}
On the other hand,
$$
\left| f\left( \frac{t}{B_{n+1}}\right) \right|
^{k(n)}=|f_{n+1}(t)|^{\rho _n}\ ,
$$
so that (10) and (13) imply (11).

Given $\varepsilon >0$, we find, for any fixed $t\in Q_p$, such
$n_0$ that
$$
|f_n(\gamma _nt)-g(\gamma _nt)|<\varepsilon \ \ \mbox{if }\
n\ge n_0
$$
(since the sequence $\{\gamma _nt\}_{n\ge 0}$ is pre-compact).
Thus
$f_n(\gamma _nt)\to g(\gamma _0t)$ by continuity of $g$, and
(12) follows from (11).  $\Box $

\begin{cor}
If $\beta \ne 0$ then $g(t)\ne 0$ for any $t\in Q_p$.
\end{cor}

{\it Proof.} Suppose that $g(t_0)=0$ for some $t_0\in Q_p$. By
(12) we find that
$$
|g(t_0)|=|g(\gamma _0t_0)|^{1/\beta }=|g(\gamma
_0^2t_0)|^{1/\beta ^2}=\ldots =|g(\gamma _0^nt_0)|^{1/\beta
^n}.
$$
Since $\gamma _0^n\to 0$ and $g(0)=1$, we obtain that
$g(\gamma _0^nt_0)\ne 0$ for a certain $n$, so that we come to
a contradiction. $\Box $

\section{Distributions of Stable Type}

In this section we consider certain distributions on $Q_p$ seen
as counterparts of classical stable distributions.

\begin{teo}
(i) Let $\Phi $ be a Borel measure on $Q_p\setminus \{0\}$
which is finite outside any neighbourhood of zero and satisfies
the relation
\begin{equation}
\Phi (M)=\beta \Phi (\gamma _0M)\ ,
\end{equation}
with $0<\beta <1,\ \gamma _0\in Q_p\ ,\ 0\ne |\gamma _0|_p\le
p^{-1}$, for any compact open subset $M\subset Q_p\setminus
\{0\}$. Then a function $g(t)$ of the form
\begin{equation}
g(t)=\exp \int \limits _{ Q_p\setminus \{0\}}(\chi (ty)-1)\Phi
(dy)
\end{equation}
is a characteristic function of a distribution which is a weak
limit of some sequence (1) with $\rho _n\to \beta ,\ \gamma
_n\to \gamma _0$.

(ii) If distributions $F_n$ for a sequence (1) with independent
symmetric identically distributed random variables $X_n$,
$|B_n|_p\to \infty ,\ \gamma _n\to \gamma _0$, $ 0\ne
|\gamma _0|_p\le p^{-1}$, $\rho _n\to \beta ,\ 0<\beta <1$,
converge weakly to a distribution $G$. then its characteristic
function is of the form (15) where the L\'evy measure $\Phi $
is symmetric and satisfies (14).
\end{teo}

{\it Proof. } (i) By [PRV], the function (15) is a
characteristic function corresponding to a random variable $X$.
Let $X_1,X_2,\ldots $ be independent copies of $X$. Set
$B_n=\gamma _0^{-n}$, $k(n)=\left[ \beta ^{-n}\right] $ where
$[\cdot ]$ means the integer part. Then
$$
f_n(t)=\left( g(\left( \gamma _0^nt\right) \right) ^{\left[
\beta ^{-n}\right]}
$$
$$
=\exp \beta ^n\left[ \beta ^{-n}\right]\int \limits _{ Q_p\setminus
\{0\}}(\chi (ty)-1)\Phi (dy)\longrightarrow g(t)\ ,\ \ n\to
\infty ,
$$
uniformly on compact subsets, since $\left[ \beta ^{-n}\right]
=\beta ^{-n}+O(1),\ n\to \infty $.

(ii) Let us proceed from the relation (10). Since $f$ is
real-valued, continuous, and $f(0)=1$, we see that $f_n(t)>0$
for each fixed $t$, if $n$ is large enough. Hence, $g(t)\ge 0$,
and by Corollary 1, $g(t)>0$ for all $t\in Q_p$. The sequence
$\{\log f_n(t)\}$ is bounded, uniformly with respect to $t$
from any compact subset of $Q_p$. Thus
$$
k(n)\left( \left( f_n(t)\right) ^{1/k(n)}-1\right) =k(n)\left(
\exp \left( \frac{1}{k(n)}\log f_n(t)\right) -1\right)
$$
$$
=\log f_n(t)+O\left( \frac{1}{k(n)}\right)\ ,\ \ n\to \infty ,
$$
so that
\begin{equation}
k(n)\left( \left( f_n(t)\right) ^{1/k(n)}-1\right)
\longrightarrow \log g(t)\ ,\ \ n\to \infty ,
\end{equation}
uniformly on compact subsets of $Q_p$.

Introducing the measures
\begin{equation}
\Phi _n(dy)=k(n)F(d(B_ny))
\end{equation}
where $F$ is the distribution of each random variable
$X_1,\ldots ,X_n,\ldots $, we may rewrite (16) as
\begin{equation}
\int \limits _{ Q_p\setminus \{0\}}(\chi (ty)-1)\Phi _n(dy)
\longrightarrow \log g(t).
\end{equation}
Note that for each $t\ne 0$ the integral in the left-hand side
of (18) is actually taken over the set of those $y$ for which
$|y|_p>|t|_p^{-1}$.

Denote the left-hand side of (18) by $\varphi _n(t)$. If $M$ is
a compact open subset of $Q_p\setminus \{0\}$, and $m$ is an
inverse Fourier transform of the indicator $\omega _M$, then by
Lemma 2
\begin{equation}
\Phi _n(M)=\int \limits _{Q_p}\varphi _n(t)m(t)\,dt\ ,
\end{equation}
and by virtue of (18)
\begin{equation}
\int \limits _{ Q_p\setminus \{0\}}\omega _M(x)\Phi
_n(dx)\longrightarrow \int \limits _{Q_p}m(t)\log g(t)\,dt\ ,\
\ n\to \infty .
\end{equation}

It follows from (20) that the sequence
$$
\int \limits _{ Q_p\setminus \{0\}}\omega (x)\Phi _n(dx)
$$
converges for any locally constant function $\omega $ with a
compact support not containing the origin. Every continuous
function with a compact support on $Q_p\setminus \{0\}$ can be
approximated uniformly by such functions (see [VVZ]). By (19),
the sequence of measures $\{\Phi _n\}$ is bounded on compact
subsets of $Q_p\setminus \{0\}$. This means that  $\{\Phi _n\}$
is a Cauchy sequence with respect to the vague topology [He],
which is sequentially complete. Thus $\Phi _n$ is vaguely
convergent to a symmetric Radon measure $\Phi $ on
$Q_p\setminus \{0\}$.

Now, in order to prove the representation (15), it is
sufficient to show that $\Phi _n\to \Phi $ in the weak sense on
each set $M_{i,\infty }=\left\{ x\in Q_p\ :\ |x|_p>p^i\right\}\
,\ i\in {\bf Z}$. By Theorem 1.1.9 of [He], that will be proved
if we show that
\begin{equation}
\Phi _n(M_{i,\infty })\longrightarrow \Phi (M_{i,\infty
})<\infty\ ,\ \ n\to \infty.
\end{equation}
Simultaneously (21) would imply the required finiteness of
$\Phi $ outside any neighbourhood of the origin.

Consider the set
$$
M_{l,i}=\left\{ x\in Q_p\ :\ p^{i+1}\le |x|_p\le p^l\right\}\
,\ \ l>i.
$$
Let us compute $\Phi _n(M_{i,l})$ using (19) where $m(t)$
corresponds to the set $M_{i,l}$. This set is a set-theoretic
difference of two balls. The Fourier transform of the indicator
function of a ball is computed in [VVZ]. Thus
$m(t)=m_l(t)-m_i(t)$ where
$$
m_j(t)=\left\{ \begin{array}{rl}
1, & \mbox{if }|t|_p\le p^{-j} \\
0, & \mbox{if }|t|_p>p^{-j}
\end{array} \right.
$$
and we obtain that
\begin{equation}
\Phi _n(M_{i,l})=-\int \limits _{p^{-l+1}\le |t|_p\le
p^{-i}}\varphi _n(t)\,dt
\end{equation}
whence
\begin{equation}
\Phi _n(M_{i,\infty })=-\int \limits _{|t|_p\le
p^{-i}}\varphi _n(t)\,dt
\end{equation}

It follows from (18) and (22) that
$$
\Phi (M_{i,l})=-\int \limits _{p^{-l+1}\le |t|_p\le
p^{-i}}\log g(t)\,dt\ .
$$
This yields
$$
\Phi (M_{i,\infty })=-\int \limits _{|t|_p\le
p^{-i}}\log g(t)\,dt\ .
$$
Comparing with (23) we come to (21).

It remains to prove the relation (14). As we have seen,
$\varphi _n(t)\sim \log f_n(t)$, $n\to \infty $, uniformly
on compact subsets of $Q_p$. Thus by Proposition 1,
\begin{equation}
\varphi _n(\gamma _nt)\longrightarrow \beta \log g(t)\ ,\ \
n\to \infty ,
\end{equation}
uniformly on compact subsets.

Let $M$ be a compact open subset of $Q_p\setminus \{0\}$. Then
\begin{equation}
\Phi (\gamma _0^{-1}M)=\lim \limits _{n\to \infty }\Phi
_n(\gamma _n^{-1}M).
\end{equation}

Indeed, let $\omega _n$ be an indicator of the set $\gamma
_n^{-1}M$, $n=0,1,2,\ldots $. Writing the action of a measure
as a functional we get
\begin{equation}
\langle \Phi _n\ ,\ \omega _n\rangle -\langle \Phi \ ,\
\omega _0\rangle=\langle \Phi _n-\Phi \ ,\ \omega _n\rangle +
\langle \Phi \ ,\ \omega _n-\omega _0\rangle .
\end{equation}

For large $n \ \ \gamma _n^{-1}M\subset M_{l',l''}$ where $l',\
l''$ are certain fixed numbers. As above, this means that the
supports of the inverse Fourier transforms $\widetilde{\omega }
_n$ of all $\omega _n$ lie in a certain compact set $N$, so that
$$
|\langle \Phi _n-\Phi \ ,\ \omega _n\rangle |\le \int \limits
_N|\varphi _n(t)-\log g(t)|\,dt\ \longrightarrow 0
$$
due to the uniform convergence. The second summand  on the right
in (26) tends to zero due to the dominated convergence theorem,
so (25) has been proved.

Next, by (19)
$$
\Phi _n(\gamma _n^{-1}M)=\int \limits _{Q_p}\varphi
_n(t)\widetilde{\omega }_n(t)\,dt\ ,
$$
$$
\widetilde{\omega }_n(t)=\int \limits _{\gamma _n^{-1}M}\chi
(-ty)\,dy=|\gamma _n|^{-1}m(\gamma _n^{-1}t)\ ,
$$
whence
$$
\Phi _n(\gamma _n^{-1}M)=\int \limits _{Q_p}\varphi _n(\gamma
_nt)m(t)\,dt\ .
$$
Now it follows from (24), (25) and (20) that
$$
\Phi (\gamma _0^{-1}M)=\beta \int \limits _{Q_p}m(t)\log g(t)\,dt
=\beta \Phi (M)
$$
which is equivalent to (14).\ \ $\Box $

\vspace{\baselineskip}
{\it Example. }\ If $|\gamma _0|_p=p^{-1}$ then the relation (14)
means that $\Phi $ is determined by its restriction to the group
of units $U_p=\{ x\in Q_p\ :\ |x|_p=1\ \}$.

For a particular example, let the above restriction be proportional
to the restriction of the Haar measure:
$$
\Phi (M_0)=a\frac{p^\alpha -1}{1-p^{-\alpha -1}}\int \limits _{M_0}dx
\ ,\ \ \alpha >0,\ a>0,
$$
for any open and closed subset $M_0\subset U_p$. Suppose that $\beta
=p^{-\alpha }$, $k(n)=\left[ p^{\alpha n}\right] $, $\gamma _0=p$.

If $M$ is a compact open subset of $Q_p\setminus \{ 0\}$ then it may
be written as a finite union
$$
M=\bigcup _NM\cap S_N)\ ,\ \ S_N=\{ x\in Q_p\ :\ |x|_p=p^N\ \} .
$$
In accordance with (14),
$$
\Phi (M)=\sum \limits _N\Phi (M\cap S_N)=\sum \limits _N
p^{-\alpha N}\Phi (p^NM\cap S_0)
$$
$$
=\sum \limits _Np^{-\alpha N}\int \limits _{p^NM\cap S_0}dx
=\sum \limits _Np^{-(\alpha +1)N}\int \limits _{M\cap S_N}dy\ ,
$$
so that
$$
\Phi (M)=a\frac{p^\alpha -1}{1-p^{-\alpha -1}}\int \limits
_{M}|x|_p^{-\alpha -1}\,dx\ .
$$

Now the identity
$$
|t|_p^\alpha =\frac{1-p^\alpha }{1-p^{-\alpha -1}}\int \limits
_{Q_p}|s|_p^{-\alpha -1}(\chi (st)-1)\,ds
$$
(see [VVZ]) shows that in this case the limit characteristic
function $g(t)$ coincides with the function $g_{a,\alpha }(t)$
mentioned in the introduction.

\vspace{\baselineskip}
Let us consider conditions for the weak convergence of the
sequence (1) with $|B_n|_p\to \infty $.

\begin{teo}
In order that the sequence (1) be weakly convergent, it is
sufficient that the measures (17) converge weakly on each set
$M_{i,\infty }\ ,\ i\in {\bf Z}$, to a measure $\Phi $ on
$Q_p\setminus \{ 0\}$, finite outside any neighbourhood of
zero. If the random variables $X_1,X_2,\ldots $ are symmetric
and $\beta \ne 0$, this condition is also necessary.
\end{teo}

{\it Proof. } The necessity was proved in the course of proving
Theorem 1. To prove the sufficiency, write $f_n(t)$ in the form
$$
f_n(t)=\left( 1+\frac{1}{k(n)}\int \limits _{Q_p\setminus \{
0\}}(\chi (ty)-1)\Phi _n(dy)\right) ^{k(n)}\ .
$$

For every fixed $t$ we have
$$
\int \limits _{Q_p\setminus \{ 0\}}(\chi (ty)-1)\Phi _n(dy)
\longrightarrow \int \limits _{Q_p\setminus \{
0\}}(\chi (ty)-1)\Phi (dy)
$$
since $\chi (ty)-1=0$ when $|y|_p\le |t|_p^{-1}$. Recalling
that $\left( 1+\displaystyle\frac{z}{k(n)}\right) ^{k(n)}\to e^z$ uniformly
on compact sets, we find that
\begin{equation}
f_n(t)\longrightarrow \exp \left( \int \limits _{Q_p\setminus \{
0\}}(\chi (ty)-1)\Phi (dy)\right)\ ,\ \ n\to \infty .
\end{equation}
Since the function in the right-hand side of (27) is
continuous, Theorem 3.3.1 of [G] implies weak convergence of
$S_n$. \ \ $\Box $

{\it Remark. } Theorems 1 and 2 can be extended easily to the
case of general non-Archimedean local fields.

\section{Degenerate Cases}

Let us consider the ``extreme'' cases $\beta =0$ and $\beta =1$
in the weakly convergent scheme (1).

\begin{prop}
If $\beta =1$ then $G$ is degenerate. If $\beta =0$ then either
$G$ is degenerate or its characteristic function $g$ has a
compact support. If  $\beta =0$ and $|\gamma _0|_p=p^{-1}$ then
$g(t)$ coincides with the right-hand side of (5).
\end{prop}

{\it Proof. } Let $\beta =1$. As before, we may assume that
$\gamma _n \to \gamma _0$, $|\gamma _0|_p\le p^{-1}$. By
Proposition 1, we have
$$
|g(\gamma _0t)|=|g(t)|
$$
for any $t\in Q_p$, so that
$$
|g(t)|=|g(\gamma _0^nt)|\ ,\ \ n=1,2,\ldots .
$$
Since $\gamma _0^n\to 0$ for $n\to \infty $, $g$ is continuous
and $g(0)=1$, this implies the identity $|g(t)|\equiv 1$.
Hence, by Lemma 1 $G$ is degenerate.

If $\beta =0$ then the same reasoning shows that $|g(\gamma
_0t)|=1$ as soon as $g(t)\ne 0$. Thus if $g(t)\ne 0$ for all
$t$ then $G$ is degenerate. Otherwise $g$ has a compact support
in accordance with (6). The last assertion of the proposition
follows from Lemma 1.\ \ $\Box $

Finally, let us consider the case of a weakly convergent
sequence (1) with $|B_n|_p\not\to \infty $. Then there exists a
subsequence $B_{n_l}\to b\in Q_p$.

\begin{prop}
If $b\ne 0$ then $g$ coincides with the right-hand side of (5).
If $b=0$ then $G$ is degenerate.
\end{prop}

{\it Proof. } Let  $b\ne 0$. It follows from (10) that the
inequality $|f(b^{-1}t)|<1$ implies the equality $g(t)=0$. If
$|f(b^{-1}t)|=1$ for some $t$ then Lemma 1 shows that $f$ is
locally constant. Thus $f(s)|=1$ when $s$ belongs to a certain
neighbourhood of the point $b^{-1}t$; in particular,
$|f(B_{n_l}^{-1}t|=1$ for large $l$ whence $|g(t)|=1$. It
remains to use Lemma 1.

Let $b=0$. Take, for an arbitrary $s\in Q_p$, a compact set $C$
containing the subsequence $\{ B_{n_l}s\}$. Given $\varepsilon
>0$, we find such a natural number $l_0$ that for $l\ge l_0$
$$
\left| \ \left| f\left( \frac{t}{ B_{n_l}}\right)
\right|^{n_l}-|g(t)|\ \right| <\varepsilon \ ,\ \ t\in C.
$$
In particular, for $t= B_{n_l}s$ we obtain that
$$
\left| \ |f(s)|^{n_l}-|g(B_{n_l}s|\ \right| <\varepsilon \ ,\ \
l\ge l_0\ .
$$
Since $g(B_{n_l}s)\to 1$, we see that $|f(s)|^{n_l}\to 1$ for
$l\to \infty $, whence $|f(s)|\equiv 1$. By (10) $G$ is
degenerate.\ \ $\Box $

\end{document}